\documentclass[12pt,leqno]{amsart}
\usepackage {enumerate}
\usepackage{amssymb}
\usepackage{amsmath}
\usepackage{amsthm}
\usepackage{pstricks, pst-poly}
\usepackage{amsbsy}
\usepackage{bm}
\usepackage{graphicx}%
\usepackage{color}
\usepackage{hyperref}
\usepackage{float}
\usepackage{appendix}
\usepackage{epstopdf}

\theoremstyle{plain}

\newtheorem{open problem}{Open Problem}

\theoremstyle{definition}

\theoremstyle{remark}

\newtheorem*{acknowledgements}{Acknowledgements}

\def\@enum@{\list{\csname label\@enumctr\endcsname}%
           {\usecounter{\@enumctr}\def\makelabel##1{\hss\llap{##1}}
           \itemsep=2pt\parsep=0pt\topsep=3pt plus 1pt minus 1 pt}}

\newenvironment{proclaim}[1]{
                \par\vspace{\topsep}\noindent{\bf #1}
                \begin{em}}
                {\end{em}\par\vspace{\topsep}}

 \newenvironment{proclama not emphasized}[1]{\par\vspace{\topsep}\noindent{\bf #1}}{\par\vspace{\topsep}}

\newcommand{\C}{\mathbb{C}}

\newcommand{\tr}{\mathrm{tr}}

\begin{document} 

\title{Self-intersection numbers of length-equivalent curves on surfaces}

\author{Moira Chas}

\address{	Department of Mathematics,\\
	Stony Brook University\\
	Stony Brook, NY, 11794}

\email{moira@math.sunysb.edu}

\date{\today}


\thanks{Partially supported by NSF grant  1098079-1-58949}

\begin{abstract}  Two free homotopy classes of closed curves in an orientable surface with negative Euler characteristic are said to be length equivalent if for any hyperbolic structure on the surface, the length of the geodesic in one class is equal to the length of the geodesic  in the other class.  We show that there are elements in the free group of two generators 
that are length equivalent and have different self-intersection numbers as elements in the fundamental group of the  punctured torus and as elements in the  pair of pants. This result answers open questions about length equivalence classes and raises new ones.

\end{abstract}

\maketitle

Consider an orientable surface $S$ (with or without boundary)with negative Euler characteristic.
A free homotopy class of curves  on a  $S$ corresponds to a conjugacy class in the fundamental group of $S$. If $S$ is endowed with a complete hyperbolic metric $m$ with geodesic boundary, each free homotopy class $x$ gets assigned a positive real number $m(x)$, the length of the unique geodesic representative in $x$  (with respect to $m$). A free homotopy class has a \emph{self-intersection number}, that is, the smallest number of crossings of representatives in general position (here, general position means that all intersection points are transversal double points). 

Two free homotopy classes $x$ and $y$ are \emph{length equivalent} if for every hyperbolic metric on  $S$, the length of the geodesic representative in $x$ equals the length of the geodesic representative in $y$. 

Two elements $X$ and $Y$ in $\pi_1(S)$ are \emph{trace equivalent} if for any representation of $\pi_1(S)$ into $SL(2,\C)$, the images of  $X$ and $Y$ have the same trace squared.

Leininger \cite[Proposition 3.2]{lein} showed that length-equivalence and trace-equivalence define the same relations. 

This note addresses the relation between self-intersection and length equivalence, by verifying following result.

\begin{proclaim}{Theorem} There exist  elements in the free group on two generators (see Table~\ref{counterex}) which are length equivalent (that is, they have the same trace squared for any representation  of the group into $SL(2,\C)$) and have different self-intersection numbers as closed curves    on the punctured torus and and as closed curves on the pair of pants.
\end{proclaim}

\begin{table}[htdp]
\begin{center}
\begin{tabular}{|c|c|c|}
\hline
Cyclically reduced  & Self-intersection Number  & Self-intersection Number \\
 word&  Pair of Pants & Punctured Torus\\
\hline
 $aaabaaBAbAABabaB$ &15& 34\\
 $aaabaBaabaBAAbAB$ &19 &32\\
 \hline
\end{tabular}
\end{center}
\caption{Length equivalent elements with different self-intersection  numbers (capital letters are used to represent inverses)}
\label{counterex}
\end{table}%

Horowitz \cite{hor}, answering a question of  Magnus, proved that in any free group of rank at least two, there exist arbitrarily large subsets of elements which are not conjugate and yet,  have the same trace for every representation of the group in  $SL(2,\C)$. 
Using Horowitz's results, Randol \cite{ran} showed that  for each positive integer $n$, there are length equivalence classes containing at least $2^n$ elements. 

Randol's result is surprising. The mystery is a bit resolved when one studies one of the several   algorithms to find length equivalence classes (see \cite{and} for a survey on this topic). These algorithms use   basic facts about traces in $SL(2,\C)$ to construct such elements, namely for each pair of matrices $A, B \in SL(2,\C)$
$$
\tr(AB)+\tr(AB^{-1})=\tr (A) \tr (B) $$$$
\tr(BAB^{-1})=\tr (A) $$$$ \tr( \mathrm{Id}) =2.
$$

Fricke \cite{fri} and Vogt \cite{vog} (see also \cite{hor}) proved that for any element $w$ in a free group of finite rank, there exists a polynomial in several variables such that the trace of $w$ under any representation of the group into $SL(2,\C)$ is equal to that polynomial evaluated at traces of  certain  products of the images of the generators.  

These trace equivalence classes(i.e., length equivalence classes) are still not completely understood  \cite{lein}. There is no known characterization from a geometric point of view. Note that the definition of trace equivalence is completely algebraic and does not distinguish different surfaces with boundary with the same Euler characteristic.

Hamenst\"adt asked at the Workshop on Kleinian Groups and Hyperbolic 3-Manifolds, held at the University of Warwick in September 2001, whether there a connection between the size of the length equivalence classes and the self-intersection numbers of the elements \cite{and}.
 
Humphries conjectured \cite[Conjecture 1.5.a]{hum} that if two elements in a free group are trace-equivalent, then they have the same self-intersection number. Note that our results is a counterexample to this conjecture.
 
Horowitz gave an algorithm that in  step $n$ yields $2^n$ non-conjugate, length-equivalent elements in the free group on two generators \cite[Example 8.2]{hor} .
We tested  the self-intersection numbers of the elements  constructed by the Horowitz algorithm for $n \in \{1,2,3,4,5,6\}$.  We computed  that all the length equivalent elements constructed by Horowitz have the same self-intersection number for $n \in \{1,2,3,4,5,6\}$ (see Table \ref{table hor})

Buser \cite[Section 3.7]{bus} gave another algorithm for finding non-conjugate, length equivalent elements in the free group on two generators.  All the length equivalent elements constructed by Buser have the same self-intersection number for $n \in \{1,2,3,4,5,6,7\}$
(see Table~\ref{table buser}.)

Computations with Horowitz's and Buser's algorithms lead us to the conjecture below, which is the focus of ongoing work with Daniel Levine and Shalin Parekh.

\begin{proclaim}{Conjecture} The $2^n$ trace-equivalent elements of step $n$ in Horowitz' algorithm have the same self-intersection number.
The $2^n$ trace-equivalent elements of step $n$ in Buser' algorithm have the same self-intersection number.
\end{proclaim}

But, because of the example in Table~\ref{counterex} there must be other methods besides these for generating length equivalence classes. 

\begin{proclaim}{Problem: } Find algorithms to generate complete length equivalence classes.
\end{proclaim}

The two classes of Table~\ref{counterex} were found with the help of a computer. First, the pair of pants was given a generic hyperbolic metric. Then  the set of all cyclically reduced words of a given word length,  was  divided into subsets that had  geometric length close enough for the chosen metric (since one needs to approximate to perform these computations, and cannot require the length to be equal but only ''close enough"). Among those subsets, the ones containing  classes with different self-intersections in the torus and in the pair of pants were chosen. 
Then those subsets of classes that were close enough in one metric, and have different self-intersection number,  were tested with a different metric, dividing them into subsets of words with  length "close enough'' in both metrics.  Next, this subclasses were divided into  sub-subclasses with the same Fricke polynomial.  Hence, the examples of Table~\ref{counterex}.
 
The pair of pants is obtained by labeling alternating edges of   an octagon by the letters $a, A, b, B$ (capital letters are used to represent inverses), and identifying edges with the same letter (without creating M\"obius bands). The generator $a$ in the pair of pants  has a representative that crosses the edge labeled  $a$ (and ''reenters" the pants through the edge labeled  $A$). Analogously, the generator $b$ has a representative that crosses the edge labeled  $b$. There is a bijection between cyclically reduced words on the $\{a,b,A,B\}$ alphabet and  a conjugacy classes in the fundamental group of the pair pants  obtained after identifying appropriate edges.

Similarly, the torus with one boundary component is obtained  by labeling alternating edges of   an octagon by the letters $a, b, A, B$. The generators of the punctured torus fundamental group are determined in the same way as those for the pair of pants. There bijection between cyclic reduced words on the $\{a,b,A,B\}$  alphabet and a conjugacy classes of the fundamental group of the torus with one boundary component so obtained. 
 
It is not hard to see that if one labels the octagon yielding the torus with the letters  $A, b, a, B$ instead of  $a, b, A, B$, the self-intersection number of a word will give the same in both cases (because switching $a$ and $A$, or $b$ by $B$ is equivalent to  performing a symmetry  on the punctured torus).   On the other hand, in the pair of pants, switching $a$ and $A$ may change self-intersection (for instance, consider the words $ab$ and $aB$; one has a simple representative, the other has a representative that is a figure eight). However, in the words we found, switching $a$ and $A$ or $b$ and $B$ does not alter the self-intersection.

Two free homotopy classes of curves $x$ and $y$ on a surface $S$ are \emph{simple-intersection equivalent} if for any simple free homotopy class $s$ in $S$, the intersection of $x$ and $s$ is equal to the intersection of $y$ and $s$. Leininger \cite[Theorem1.4]{lein} proved that length-equivalence implies simple intersection equivalence. Thus our examples imply the following corollary, (compare Leininger \cite{lein})

\begin{proclaim}{Corollary} There are free homotopy classes of curves which are simple-intersection equivalent and have different self-intersection number.
\end{proclaim}

\begin{acknowledgements} The classes exhibited in this paper were found using a parameterization of the pair of pants we learned from Bernie Maskit. We use a Mathematica program of Goldman to compute Fricke polynomials. Cameron Crowe provided very valuable help in modifying Goldman's program to suit the needs of this study. Kaiqiao Li helped us to program the Horowitz algorithm and Daniel Levine helped  to the program Buser's algorithm. Dennis Sullivan and Anthony Phillips gave us valuable comments for this manuscript.
\end{acknowledgements}

\appendix
\section{Computations}
\subsection{Intersection numbers}
The intersection of the cyclically reduced words  below can be computed with Cohen-Lustig \cite{cl} or Arettines \cite{aret} algorithms.

\begin{table}[htdp]
\begin{center}
\begin{tabular}{|c|c|c|c|}
\hline
n  & Self-intersection Number  in & Self-intersection Number in & Word length\\
 word& Punctured Torus & Pair of Pants& \\
\hline
1& 5&4&7\\
2&83&47&29\\
3&1301&725&111\\
4&20759&11543&433\\
5&332057&184601&1715\\
6&5312795&2953499&6837\\
 \hline
\end{tabular}
\end{center}
\caption{Self-intersection of the length equivalent elements in Horowitz  algorithm}
\label{table hor}
\end{table}%

\begin{table}[htdp]
\begin{center}
\begin{tabular}{|c|c|c|c|}
\hline
n  & Self-intersection Number  in & Self-intersection Number in & Word length\\
 word& Punctured Torus & Pair of Pants& \\
\hline
1& 4 & 1&5\\
2&48&10&17\\
3&476&91&53\\
4&4432&820&161\\
5&40356&7381&485\\
6& 364640& 66430& 1457\\
7& 3286108&597871 &4373 \\
 \hline
\end{tabular}
\end{center}
\caption{Self-intersection of the length equivalent elements in Buser algorithm}
\label{table buser}\end{table}%

\subsection{Fricke Polynomial}
The Fricke polynomial of both $$aaabaaBAbAABabaB \mbox{ and } aaabaBaabaBAAbAB$$  is the following
$$
-x^8 y^2 z^2+x^7 y^3 z^3+2 x^7 y^3 z+2 x^7 y z^3-x^7 y z-3 x^6 y^4 z^2-x^6 y^4-3 x^6 y^2 z^4+4 x^6 y^2 z^2+x^6 y^2$$$$
-x^6 z^4+x^6 z^2+3 x^5 y^5 z+5 x^5 y^3 z^3-12 x^5 y^3 z+3 x^5 y z^5-12 x^5 y z^3+5 x^5 y z-x^4 y^6+6 x^4 y^4+$$$$
7 x^4 y^2 z^2-5 x^4 y^2-x^4 z^6+6 x^4 z^4-5 x^4 z^2-3 x^3 y^5 z-6 x^3 y^3 z^3+10 x^3 y^3 z-3 x^3 y z^5+10 x^3 y z^3-3 x^3 y z$$$$
+x^2 y^6+3 x^2 y^4 z^2-5 x^2 y^4+3 x^2 y^2 z^4-10 x^2 y^2 z^2+3 x^2 y^2+x^2 z^6-5 x^2 z^4+3 x^2 z^2+x^2-x y z+y^2+z^2-2
$$

\section{Representatives of the classes}
\begin{figure}[htbp]
 \begin{center}
\includegraphics[scale=0.5]{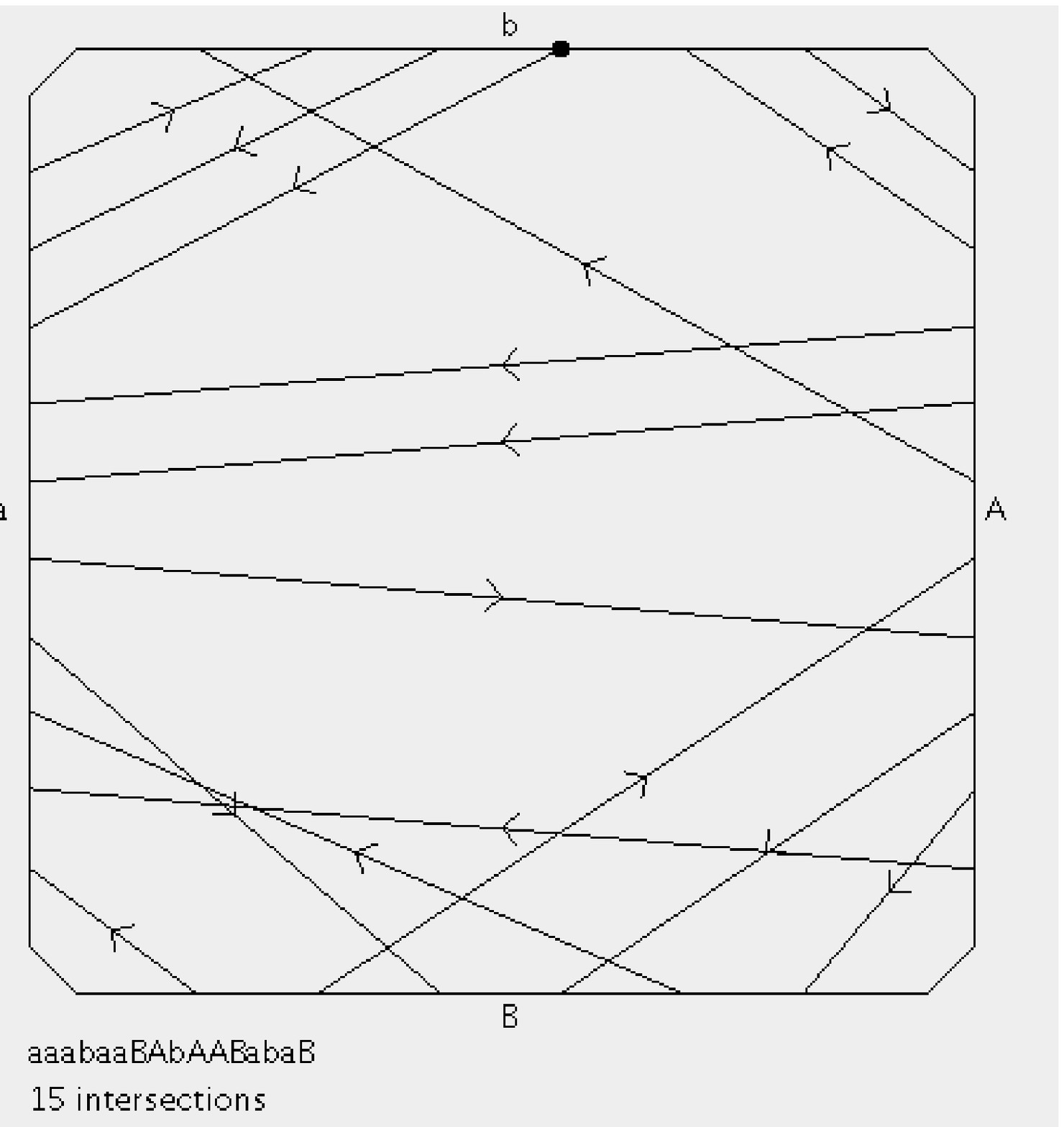}
\includegraphics[scale=0.5]{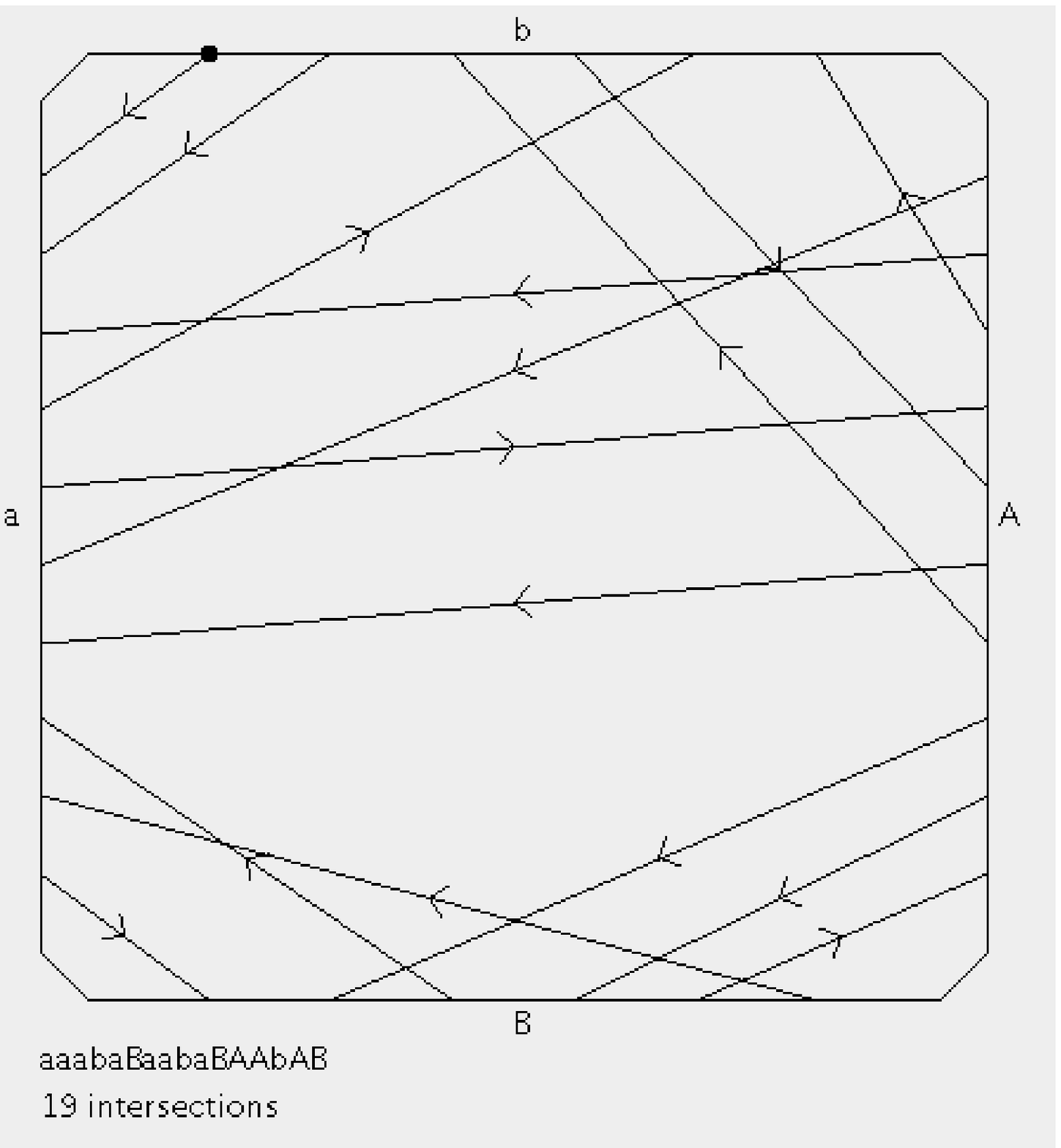}
\includegraphics[scale=0.5]{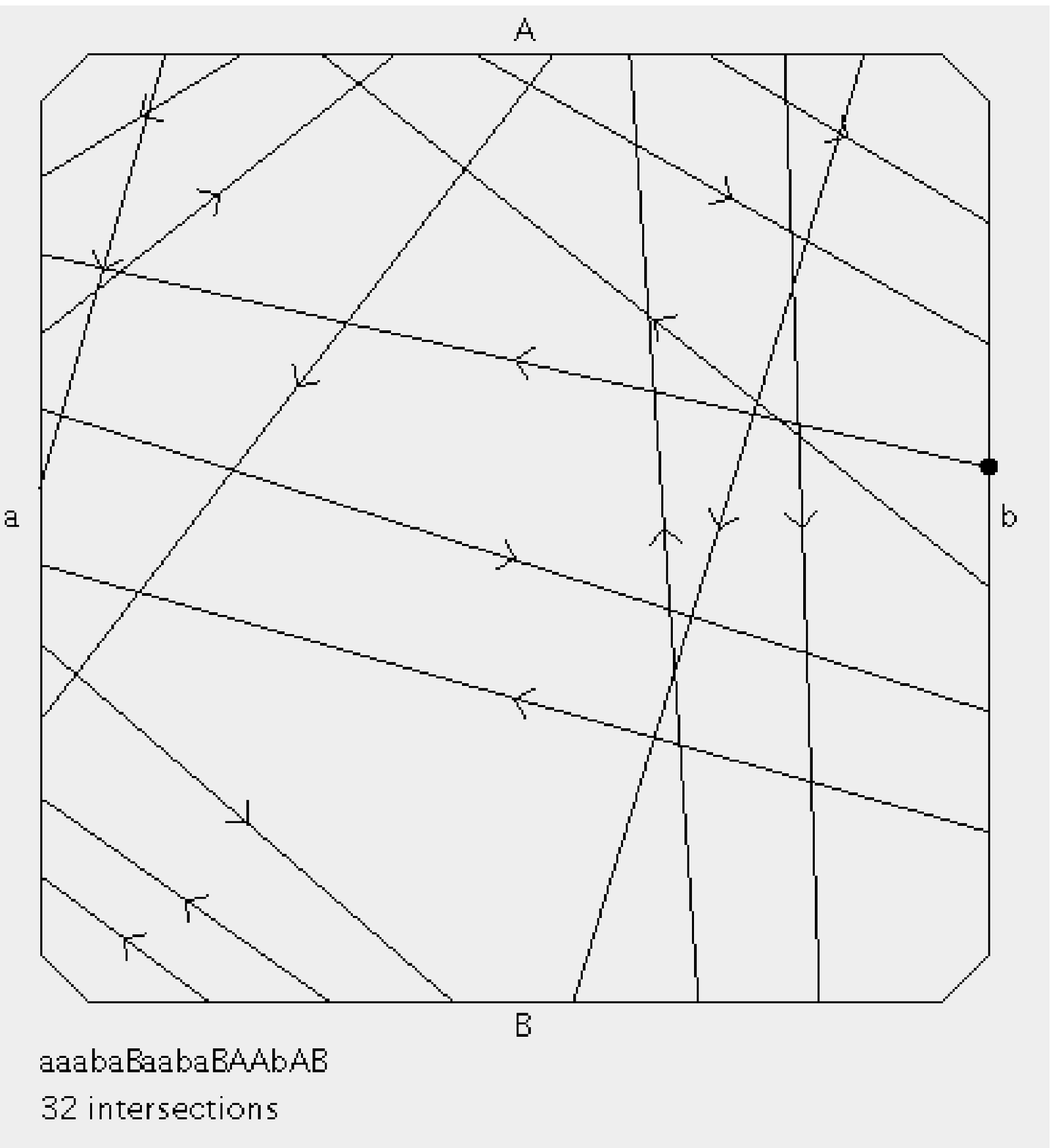}
\includegraphics[scale=0.5]{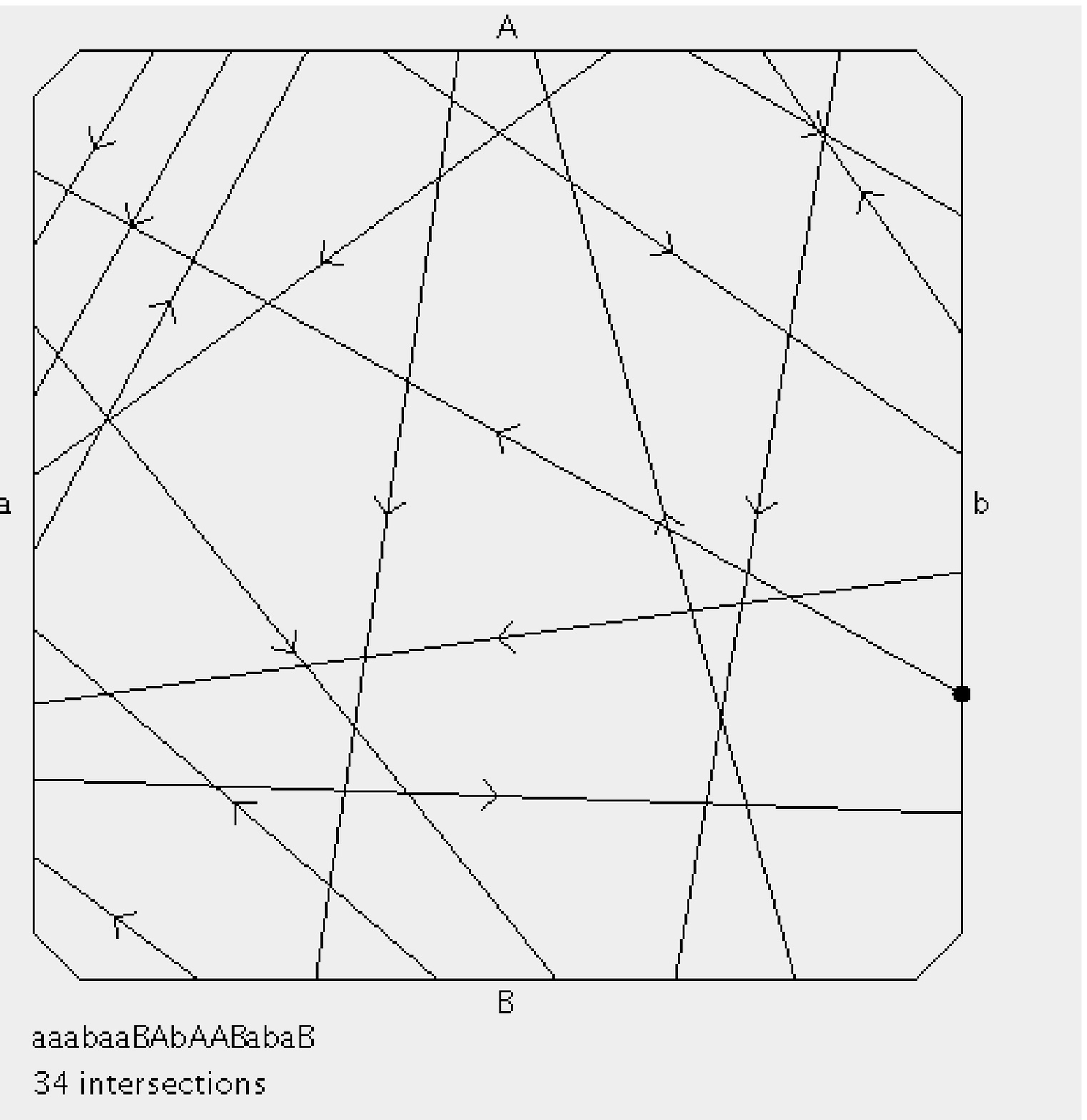}
 
   \caption{Representatives of the examples}
   \label{fh}
   \end{center}
 \end{figure}

\enddocument